\documentclass[12pt]{amsart}
\usepackage{amsmath}
\usepackage{amssymb}
\usepackage{stmaryrd}
\usepackage{epsfig,color,esint}

\numberwithin{equation}{section}

\newtheorem{prop}{Proposition}[section]
\newtheorem{theo}{Theorem}[section]
\newtheorem{lemm}{Lemma}[section]

\def\begeq{\begin{equation}}
\def\endeq{\end{equation}}

\begin{document}

\title{Li-Yau Inequality and Liouville property to a Semilinear Heat Equation on Riemannian Manifolds}
\author{Huan-Jie Chen$^1$, Shi-Zhong Du$^2$ and Yue-Xiao Ma$^3$}
\thanks{The authors are partially supported by NSFC (12171299), and GDNSF (2019A1515010605)}

    \address{The Department of Mathematics,
            Shantou University, Shantou, 515063, P. R. China.} \email{19hjchen@stu.edu.cn}

  \address{The Department of Mathematics,
            Shantou University, Shantou, 515063, P. R. China.} \email{szdu@stu.edu.cn}

  \address{The Department of Mathematics,
            Shantou University, Shantou, 515063, P. R. China.} \email{21yxma@stu.edu.cn}

\renewcommand{\subjclassname}{%
  \textup{2010} Mathematics Subject Classification}
\date{May. 2023}

\begin{abstract}
  Entire solutions of the nonlinear equation
     \begin{equation}\label{e0.1}
       u_t-\triangle_g u=u^p, \ \ x\in M, t\in{\mathbb{R}}
     \end{equation}
   were studied, where $M$ is a $N$ dimensional complete Riemannian manifold equipped with the metric $g$ and $p$ is assumed to be greater than one. The first part of this paper is devoted to investigation of Liouville property of
      \begin{equation}\label{e0.2}
        u_t-\triangle_gu=f(u), \ \ x\in M, t\in{\mathbb{R}}
      \end{equation}
  on compact manifolds, which extends a result by Castorina-Mantegazza \cite{CM2} for positive $f$. Secondly,  we will turn to non-compact manifolds and prove a Liouville theorem of \eqref{e0.1} under the assumptions of boundedness of the Ricci curvature from below, diffeomorphism of $M$ with ${\mathbb{R}}^N$  and sub-criticality of $p$ defined below. Finally, we also present simplified proofs of Yau's theorem for harmonic function and Gidas-Spruck's theorem for elliptic semilinear equation. Our proofs are based on Li-Yau type estimation for nonlinear equations.
 \end{abstract}

\maketitle\markboth{Li-Yau Inequality}{Liouville Theorem}

\section{Introduction}

One of the most important models in partial differential equation is the Laplace equation
  \begin{equation}\label{e1.1}
     \triangle_gu=0,
  \end{equation}
on a Riemannian manifold $M$ equipped with the metric $g$. The solutions of \eqref{e1.1} are called to be harmonic functions as usually. A first celebrated Liouville theorem for nonnegative harmonic function was proven by Yau \cite{Y} in 1975, which asserts that for complete noncompact Riemannian manifold with nonnegative Ricci curvature, every nonnegative harmonic function must be a constant. A natural question is to ask whether all nonnegative solutions to linear heat equation
   \begin{equation}\label{e1.2}
     u_t-\triangle_gu=0
   \end{equation}
are identical to constants? This may not be true in general without some additional assumptions. Actually, even for $M={\mathbb{R}}^N$, there exist counter examples defined by
   $$
     u(x,t)=e^{t+a\cdot x},\ \ a=(a_1,a_2,\cdots, a_N),\ \ |a|^2=\Sigma_{i=1}^Na_i^2=1.
   $$
However, after imposing some natural growth bounds on solutions and using a local version of Li-Yau type estimation, Souplet-Zhang \cite{SZ} have shown a similar Liouville theorem parallel to that for harmonic functions.

 Turning to the nonlinear equation of the form
   \begin{equation}\label{e1.3}
     -\triangle_gu=u^p, \ \ p>1,
   \end{equation}
   Gidas-Spruck have shown a Liouville theorem in their elegant paper \cite{GS} for Euclidean space $M={\mathbb{R}}^N$ and
     $$
      1<p<p_S(N),
     $$
 where
  $$
   p_S(N)\equiv\begin{cases}
        \frac{N+2}{N-2}, & N\geq3\\
        +\infty, & N=1, 2
      \end{cases}
  $$
 stands for the critical Sobolev exponent. The sharpness of the region of $p$ is due to the non-trivial positive solution
     $$
       u(x)=\Bigg(\frac{\sqrt{N(N-2)}}{1+|x|^2}\Bigg)^{\frac{N-2}{2}}, \ \ x\in M\equiv{\mathbb{R}}^N
     $$
   for $p=p_S(N)$ and dimensions $N\geq3$. When considering the parabolic version
     \begin{equation}\label{e1.4}
       u_t-\triangle_gu=u^p, \ \ p>1,
     \end{equation}
   a long standing conjecture asserts the validity of nonexistence of positive solutions for $1<p<p_S(N)$. Even in the Euclidean case $M={\mathbb{R}}^N$, this problem remained open for a long time. A first partial result was obtained by Bidaut-V\'{e}ron \cite{BV} for
      $$
       1<p<p_B(N)\equiv\begin{cases}
          \frac{N(N+2)}{(N-1)^2}, & N\geq2,\\
          +\infty, & N=1.
       \end{cases}
      $$
   Soon after, Quittner \cite{Q} solved completely the two dimensional conjecture in 2016. See also \cite{PQ1,PQ2} for issues of radial symmetric solutions. More recently, in a celebrating paper by Quittner \cite{Q2}, the full Liouville property of \eqref{e1.4} was solved completely for Euclidean spaces and all subcritical Sobolev exponents $1<p<p_S(N)$.

  For the Riemannian manifold case, few bibliographies were known so far. Throughout this paper, letting $T$ be a real number, we will call the solution to be an ancient solution if it is defined on $M\times(-\infty,T)$ and to be an eternal solution if its life region equals to $M\times(T,+\infty)$. Furthermore, if the solution is defined on whole space-time $M\times(-\infty,+\infty)$, we will call it to be an entire solution of \eqref{e1.1}. In the first part of this paper, we will consider the compact manifold case and extend results of Cartorina-Mantegazza \cite{CM, CM2, CMS} to a more complete situation as following.\\

 \begin{theo}\label{t1.1}
 Let $(M,g)$ be a $N$ dimensional complete compact manifold equipped with the metric $g$ and consider $u$ to be a solution (may change sign) of
  \begin{equation}\label{e1.5}
    u_t-\triangle_gu=f(u), \ \ x\in M, t\in{\mathbb{R}}.
  \end{equation}

\noindent(1) If $f$ a positive convex function on $(0,+\infty)$ satisfying
    \begin{equation}\label{e1.6}
       \int^{+\infty}_1\frac{dz}{f(z)}<+\infty,
    \end{equation}
then \eqref{e1.5} has no eternal positive solution.

\noindent(2) If $f$ is a nonnegative function satisfying that $f(0)=0, f(z)\not=0, \forall z\not=0$ and
  \begin{equation}\label{e1.7}
    \int^{+\infty}_1\Bigg(\frac{1}{f(z)}+\frac{1}{f(-z)}\Bigg)dz<+\infty,
  \end{equation}
then any nonnegative entire solution $u$ must be identical to zero.

\noindent(3) If $f(u)=|u|^{p-1}u$ for some $p\in(1,p_S(N))$, then \eqref{e1.5} admits non-trivial entire solution $u$.
\end{theo}

It is remarkable that Part (2) in the above theorem has been obtained recently by Castorina-Mantegazza \cite{CM2} for $f(z)=|z|^p, p>1$. Moreover, Part (3) gives counter examples of Liouville property for change sign $f$. At the next stage, we will turn to the non-compact manifolds and prove the following Liouville theorem for nonlinear heat equation \eqref{e1.4}:\\

\begin{theo}\label{t1.2} Let $(M,g)$ be a $N$ dimensional complete noncompact manifold equipped with the metric $g$ and $u$ be a nonnegative entire solution of \eqref{e1.4}. Supposing that
    \begin{equation}\label{e1.8}
      1<p<p_*(N)
    \end{equation}
for $p_*(N)$ defined by
   \begin{equation}\label{e1.9}
    p_*(N)\equiv\begin{cases}
         \frac{N+2+\sqrt{N^2+8N}}{2(N-1)}, & N\geq2,\\
         8, & N=1,
      \end{cases}
   \end{equation}
and the Ricci curvature satisfying
   $$
    \begin{cases}
       R_{ij}\geq0, & M\ \mbox{is an arbitrary manifold,}\\
       R_{ij}\geq -Kg_{ij}, & M\ \mbox{is a suitable manifold}
    \end{cases}
   $$
for positive constant $K$, then $u$ must be identical to zero.
\end{theo}

More recently, Liouville property to nonnegative eternal solutions of \eqref{e1.4} has been obtained by Castorina-Mantegazza-Sciunzi \cite{CMS} for nonnegative Ricci curvature and exponents $1<p\leq p_F\equiv 1+\frac{2}{N}$. When consider entire solutions of \eqref{e1.4} in \cite{CCM}, the authors also mentioned that the result of Bidaut-V\'{e}ron \cite{BV} in Euclidean spaces could be generalized to Riemaniann manifold with nonnegative Ricci curvature, using a same parabolic version of Gidas-Spruck type estimation. Our main result Theorem \ref{t1.2} enables negative Ricci curvature and relies a complete different argument. The third part of this paper is devoted to the simplified proofs of Yau's theorem for harmonic functions and Gidas-Spruck's theorem for nonlinear elliptic equation \eqref{e1.3}:\\

\begin{theo}\label{t1.3} Let $(M,g)$ be a $N$ dimensional complete noncompact manifold with the metric $g$. We have the following alternatives:

 (1) under the assumption of Theorem \ref{t1.2}, there is no positive solution to \eqref{e1.3}, and

 (2) assuming only $R_{ij}\geq0$, every nonnegative solution to \eqref{e1.1} must be a constant.\\
\end{theo}

Our proofs of Theorem \ref{t1.2} and \ref{t1.3} were based on the following differential Harnack inequalities of Li-Yau type:\\

\begin{theo}\label{t1.4} Let $(M,g)$ be a $N$ dimensional complete noncompact manifold equipped with the metric $g$. Assume that $R_{ij}\geq-Kg_{ij}$ holds for some nonnegative constant $K$, then

\noindent (1) for any positive entire solutions to \eqref{e1.2} and any given $\gamma>1$, there exists a positive constant $C_{N,\gamma}$ depending only on $N,\gamma$, such that
    \begin{equation}\label{e1.10}
      \Bigg|\frac{\nabla u}{u}\Bigg|^2-\gamma(\log u)_t\leq C_{N,\gamma}K.
    \end{equation}

\noindent (2) for any positive entire solutions $u$ to \eqref{e1.4} and $p$ satisfying the subcritical condition \eqref{e1.8}, there exist positive constants $\beta, \gamma$ and $C_0$ depending only on $N, p$, such that
   \begin{equation}\label{e1.11}
     \Bigg|\frac{\nabla u}{u}\Bigg|^2-\gamma(\log u)_t+\beta u^{p-1}\leq C_0K.
   \end{equation}
\end{theo}

\vspace{10pt}

A first related version of part (1) was obtained by Li-Yau in \cite{LY} for eternal positive solution to linear heat equation, which was called ``Li-Yau" type inequality in the literatures therein. A similar form of part (2) was thereafter extended by (J.Y.) Li \cite{L} to semilinear heat equation for $1<p<\frac{N}{N-1}$. The readers may also refer to a paper \cite{CCM2} for a weaker form of Li-Yau type inequality on Riemannian manifold with nonnegative Ricci curvature.

 The contents of this paper are organized as following: In Section 2, Theorem \ref{t1.1} will be proven for different $f$. Secondly, we show a basic gradient formula to solution of \eqref{e1.3} using the Bochner type identity. With the help of a geometric lemma, we obtain a differential Harnack inequality of Li-Yau type in Section 3 and 4. As an application, we were able to prove the non-existence of positive entire solution for nonlinear parabolic equation in Section 5. Finally, simple proofs of Liouville theorems for harmonic functions of \eqref{e1.1} and positive solutions of \eqref{e1.3} will be presented in Section 6.

\vspace{40pt}

\section{Liouville theorems on compact manifolds}

At the beginning, we assume $(M,g)$ to be a complete compact manifold equipped with the metric $g$. The proof of Part (1) of Theorem \ref{t1.1} was clear. Actually, integrating \eqref{e1.5} over manifold, we get
  $$
   \frac{d}{dt}\int_Mu=\int_Mf(u)\geq f\Big(\int_Mu\Big)
  $$
by Jensen's inequality. As a result,
   $$
    \int^{\int_Mu(T)}_{\int_Mu(T_0)}\frac{dz}{f(z)}\geq T-T_0
   $$
holds for any $T>T_0$. Sending $T\to+\infty$ and using \eqref{e1.6}, we reach a contradiction and thus complete the proof of Part (1) of Theorem \ref{t1.1}. $\Box$ \\

Part (2) is consequence of the following lemma.

\begin{lemm}\label{l2.1}
  Under the assumption of Theorem \ref{t1.1}, we have
    \begin{equation}\label{e2.1}
      \varphi(t)\equiv\min_{x\in M}u(x,t)\equiv u(\xi(t),t)\equiv 0, \ \ \forall t\in{\mathbb{R}}
    \end{equation}
 holds for some $\xi(t)\in M, \forall t$.
\end{lemm}

\noindent\textbf{Proof.} Noting that $\varphi(t)$ is a Lipschitz function in $t$ and
   \begin{equation}\label{e2.2}
      \frac{d\varphi}{dt}\geq f(\varphi),\ \ \forall t
   \end{equation}
if $\varphi(t_0)>0$ at some $t_0\in{\mathbb{R}}$, it is inferred from \eqref{e2.2} that
   \begin{equation}\label{e2.3}
     \int^{\varphi(T)}_{\varphi(t_0)}\frac{d\varphi}{f(\varphi)}\geq T-t_0, \ \ \forall T>t_0.
   \end{equation}
However, sending $T\to+\infty$ and using \eqref{e1.7}, we obtain that the L.H.S. in \eqref{e2.3} is finite but the R.H.S. in \eqref{e2.3} tends to infinite. Contradiction holds. On another hand, if $\varphi(t_0)<0$ at some $t_0\in{\mathbb{R}}$, we set $\phi(t)=\varphi(t_0-t)$ and obtain that
  \begin{equation}\label{e2.4}
    \frac{d\phi}{dt}\leq-f(\phi), \ \ \forall t.
  \end{equation}
As a result,
   \begin{equation}\label{e2.5}
    \int^{\varphi(t_0-T)}_{\varphi(t_0)}\frac{d\varphi}{f(\varphi)}=\int^{\phi(T)}_{\phi(0)}\frac{d\phi}{f(\phi)}\leq -T,\ \ \forall T>0.
   \end{equation}
Sending $T\to+\infty$ in \eqref{e2.5}, we arrive at a contradiction again by \eqref{e1.7}. The proof of the lemma was done. $\Box$\\

Using Lemma \ref{l2.1} and the strong maximum principle of heat equation, Part (2) of Theorem \ref{t1.1} follows. $\Box$\\

We are now in position to show Part (3) of Theorem \ref{t1.1}. Actually, we will show that there exists non-trivial steady state of \eqref{e1.5} by variational method. Considering the variational problem which minimizes the functional
  $$
   E(u)\equiv\frac{1}{2}\int_M|\nabla u|^2
  $$
upon the family of functions defined by
  $$
   {\mathcal{A}}\equiv\Bigg\{u\in H^1(M)\Bigg|\ \int_M|u|^{p+1}=1, \ \ \int_M|u|^{p-1}u=0\Bigg\}.
  $$
We shall need the following version of Sobolev embedding theorem on manifolds. (see Gilbarg-Trudinger \cite{GT})

\begin{lemm}\label{l2.2}
   Let $(M,g)$ be a $N$ dimensional complete compact Riemannian manifold equipped with the metric $g$. For each $p>0$, there exists a positive constant $C_{p,M}$ depending only on $p$ and the manifold $M$, such that
     \begin{equation}\label{e2.6}
       \Bigg(\int_M |h|^{\frac{2N}{N-2}}\Bigg)^{\frac{N-2}{N}}\leq C_{p,M}\Bigg(\int_M|\nabla h|^2+\Big(\int_M|h|^{p}\Big)^{\frac{2}{p}}\Bigg)
     \end{equation}
   holds for any functions $h\in H^1(M)$. Moreover, each bounded set of $H^1(M)$ is precompact in $L^q(M)$ for each $1<q<\frac{2N}{N-2}$.
\end{lemm}

Now, Part (3) of Theorem \ref{t1.1} is just a consequence of the following proposition.

\begin{prop}
  Let $(M,g)$ be a $N$ dimensional complete compact Riemannian manifold equipped with the metric $g$. If $p\in(1,p_S(N))$, there exists a non-trivial smooth solution $u$ of
     \begin{equation}\label{e2.7}
       -\triangle_gu=|u|^{p-1}u, \ \ \mbox{in } M,
     \end{equation}
  the Euler-Lagrange equation of the variational problem minimizing $E(\cdot)$ upon the family ${\mathcal{A}}$.
\end{prop}

\noindent\textbf{Proof.} As usually, taking a minimizing sequence $u_j, j\in{\mathbb{N}}$ of $E(\cdot)$, we have $u_j$ is bounded in $H^1(M)$. So, it is inferred from compact embedding Lemma \ref{l2.2} that there exists a limiting function $u_\infty\in{\mathcal{A}}$ satisfying
   \begin{equation}\label{e2.8}
    \begin{cases}
      u_j\rightharpoonup u_\infty, & \mbox{ weakly in } H^1(M),\\
      u_j\rightarrow u_\infty, & \mbox{ strongly in } L^q(M), \ \forall q\in(1,p_S(N)+1).
    \end{cases}
   \end{equation}
By lower semi-continuity of weak convergence, we also have
   \begin{equation}\label{e2.9}
     \inf_{h\in{\mathcal{A}}}E(h)\leq E(u_\infty)\leq\liminf_{j\to\infty}E(u_j)=\inf_{h\in{\mathcal{A}}}E(h).
   \end{equation}
Hence, $u_\infty$ is the desired minimizer of $E(\cdot)$ upon the family ${\mathcal{A}}$, which satisfies the Euler-Lagrange equation
   \begin{equation}\label{e2.10}
     -\triangle_gu_\infty=\lambda|u_\infty|^{p-1}u_\infty+\mu|u_\infty|^{p-1}
   \end{equation}
for some multipliers $\lambda, \mu\in{\mathbb{R}}$. We claim that $\lambda>0$ and $\mu=0$. Actually, integrating \eqref{e2.10} over $M$ and using the fact
   $$
    \int_M|u_\infty|^{p-1}u_\infty=0,
   $$
we obtained that $\mu=0$. Multiplying \eqref{e2.10} by $u_\infty$ and then integrating by parts, it is also inferred that $\lambda>0$. Now, by multiplying a constant to $u_\infty$, we obtain a desired non-trivial solution of \eqref{e2.7}. The proposition was shown. $\Box$\\

\vspace{40pt}

\section{Bochner identity and basic formulas}

Let $(M,g)$ be a Riemannian manifold equipped with a complete metric $g$. Under local coordinates $x=(x^1,\cdots, x^N)$ of $M$, the standard Laplace-Beltrami operator is defined by
    $$
     \triangle_g\equiv\frac{1}{\sqrt{\det g}}\frac{\partial}{\partial x^i}\Bigg(\sqrt{\det g}g^{ij}\frac{\partial}{\partial x^j}\Bigg).
    $$
Hereafter, without particular indication, we drop the subscript $g$ of $\triangle_g$ for simplicity. Given any function $U\in C^\infty(M)$, the Bochner identity gives
   \begin{equation}\label{e3.1}
    \triangle|\nabla U|^2=2|\nabla^2U|^2+2\nabla U\nabla\triangle U+2Ric(\nabla U,\nabla U),
   \end{equation}
where
   \begin{eqnarray*}
    R_{ij}&=&\frac{\partial}{\partial x^i}\Gamma^t_{tk}-\frac{\partial}{\partial x^t}\Gamma^t_{ik}+\Gamma^s_{tk}\Gamma^t_{is}-\Gamma^s_{ik}\Gamma^t_{ts},\\
    \Gamma^k_{ij}&=&\frac{1}{2}g^{kl}\Bigg(\frac{\partial}{\partial x^i}g_{jl}+\frac{\partial}{\partial x^j}g_{il}-\frac{\partial}{\partial x^l}g_{ij}\Bigg)
   \end{eqnarray*}
are the components of Ricci curvature tensor $Ric$ and the Christoffel symbol of $M$ respectively. Our main results were asserted for complete non-compact Riemannian manifolds with nonnegative Ricci curvature. If the Ricci curvature tensor is merely assumed to be bounded from below by a negative constant $-K$, we will use the assumption of the ``suitability" of the manifold, which means that $M$ is diffeomorphism to ${\mathbb{R}}^N$. The importance of ``suitability" of manifold $M$ is that the equation \eqref{e1.4} adhere a one-parameter symmetric group
   $$
     (x,t,u)\to (k^{-1}x, k^{-2}t, k^{\frac{2}{p-1}}u)
   $$
in its global coordinates $x\in{\mathbb{R}}^N$.

Now, letting $u$ be a positive solution to \eqref{e1.4} and defining $v\equiv\log u$, we have
   \begin{equation}\label{e3.2}
     v_t-\triangle v=u^{p-1}+|\nabla v|^2.
   \end{equation}
Setting $w=|\nabla v|^2=\big|\frac{\nabla u}{u}\big|^2$ and using Bochner identity \eqref{e3.1}, we get
   \begin{eqnarray}\nonumber\label{e3.3}
     (\partial_t-\triangle)w&=&2\nabla v\nabla v_t-2\nabla v\nabla\triangle v-2\big|\nabla^2v\big|^2-2Ric(\nabla v,\nabla v)\\ \nonumber
     &=&2\nabla v\nabla\big(u^{p-1}+|\nabla v|^2\big)-2\big|\nabla^2v\big|^2-2Ric(\nabla v,\nabla v)\\
     &=&-2\big|\nabla^2v\big|^2-2Ric(\nabla v,\nabla v)+2(p-1)u^{p-1}w+2\nabla v\nabla w.
   \end{eqnarray}
Taking derivative in time, it is inferred from \eqref{e3.2} that
   \begin{equation}\label{e3.4}
     (\partial_t-\triangle)(\log u)_t=(p-1)u^{p-1}(\log u)_t+2\nabla v\nabla(\log u)_t.
   \end{equation}
So, a combination of \eqref{e3.3} with \eqref{e3.4} shows that
  \begin{eqnarray}\nonumber\label{e3.5}
   &&(\partial_t-\triangle)\Bigg\{\Big|\frac{\nabla u}{u}\Big|^2-\gamma(\log u)_t\Bigg\}=-2|\nabla^2\log u|^2+2\nabla\log u\nabla\Bigg\{\Big|\frac{\nabla u}{u}\Big|^2-\gamma(\log u)_t\Bigg\}\\
   && \ \ \ \ \ \ \ \ \ \ \ \ \ \  -2Ric(\nabla v,\nabla v)+2(p-1)u^{p-1}\Big|\frac{\nabla u}{u}\Big|^2-\gamma(p-1)u^{p-1}(\log u)_t
  \end{eqnarray}
for any constant $\gamma$. Substituting \eqref{e3.2} into \eqref{e3.5}, we obtain that
  \begin{eqnarray}\nonumber\label{e3.6}
   &&(\partial_t-\triangle)\Bigg\{\Big|\frac{\nabla u}{u}\Big|^2-\gamma(\log u)_t\Bigg\}=-2|\nabla^2\log u|^2+2\nabla\log u\nabla\Bigg\{\Big|\frac{\nabla u}{u}\Big|^2-\gamma(\log u)_t\Bigg\}\\ \nonumber
   && \ \ \ \ \ \ \  \ \ \ \ \ \ \  \ \ \ \ \ \ \  \ \ \ \ \ \ \ -2Ric(\nabla v,\nabla v)-(\gamma-2)(p-1)u^{p-1}\Big|\frac{\nabla u}{u}\Big|^2\\
   && \ \ \ \ \ \ \  \ \ \ \ \ \ \  \ \ \ \ \ \ \  \ \ \ \ \ \ \ \ \ -\gamma(p-1)u^{p-1}\triangle\log u-\gamma(p-1)u^{2p-2}.
  \end{eqnarray}
On another hand, since
   \begin{equation}\label{e3.7}
    (\partial_t-\triangle)u^{p-1}=(p-1)u^{2p-2}-(p-1)(p-2)u^{p-1}\Big|\frac{\nabla u}{u}\Big|^2,
   \end{equation}
if one sets $\rho\equiv\Big|\frac{\nabla u}{u}\Big|^2-\gamma(\log u)_t+\beta u^{p-1}$ for some constant $\beta$, there holds
   \begin{eqnarray}\nonumber\label{e3.8}
     \rho_t-\triangle\rho&=&-2|\nabla^2\log u|^2-2Ric(\nabla v,\nabla v)+2\nabla\log u\nabla\Bigg\{\Big|\frac{\nabla u}{u}\Big|^2-\gamma(\log u)_t\Bigg\}\\ \nonumber
     &&-\gamma(p-1)u^{p-1}\triangle\log u-(\gamma-\beta)(p-1)u^{2p-2}\\ \nonumber
     &&-(p-1)\Bigg\{\gamma-2+\beta(p-2)\Bigg\}u^{p-1}\Big|\frac{\nabla u}{u}\Big|^2\\ \nonumber
     &=&-2|\nabla^2\log u|^2-2Ric(\nabla v,\nabla v)+2\nabla\log u\nabla\rho-(\gamma-\beta)(p-1)u^{2p-2}\\
     &&-\gamma(p-1)u^{p-1}\triangle\log u-(p-1)\Bigg\{\gamma-2+\beta p\Bigg\}u^{p-1}\Big|\frac{\nabla u}{u}\Big|^2.
   \end{eqnarray}
Now, we have the following lemma:

\begin{lemm}\label{l3.1}
  Suppose that Ricci curvature satisfies that $R_{ij}\geq-Kg_{ij}$ for some nonnegative constant $K$. If
    \begin{equation}\label{e3.9}
      1<p<p_*(N)\equiv\begin{cases}
         \frac{N+2+\sqrt{N^2+8N}}{2(N-1)}, & N\geq2,\\
         8, & N=1,
      \end{cases}
    \end{equation}
  there exist positive constants $\beta, \gamma$ satisfying
    \begin{equation}\label{e3.10}
      \gamma>\max\{1,\beta,2-\beta p\},\ \ \ \frac{\gamma^2}{\gamma-\beta}(p-1)<\frac{8}{N}
    \end{equation}
  and positive constants $\alpha_1,\alpha_2, \alpha_3$ depending only on $N,p,\beta,\gamma$ such that
    \begin{equation}\label{e3.11}
      \rho_t-\triangle\rho\leq-\alpha_1|\triangle\log u|^2+2\nabla\log u\nabla\rho-\alpha_2\rho^2-\alpha_3\Big|\frac{\nabla u}{u}\Big|^4+2K\Big|\frac{\nabla u}{u}\Big|^2
    \end{equation}
  holds for positive solution $u$ of \eqref{e1.4} and $\rho\equiv\Big|\frac{\nabla u}{u}\Big|^2-\gamma(\log u)_t+\beta u^{p-1}$ defined as above.
\end{lemm}

\noindent\textbf{Proof.} The proof is element and we will focus on the case $N\geq2$. By Cauchy's inequality, we have
   $$
    |\nabla^2\log u|^2\geq\frac{1}{N}|\triangle\log u|^2.
   $$
  So, to show the validity of \eqref{e3.11} under \eqref{e3.9}, we need only to choose $\beta,\gamma$ satisfying \eqref{e3.10} such that
   \begin{equation}\label{b.1}
    -\frac{2}{N}|\triangle\log u|^2-(\gamma-\beta)(p-1)u^{2p-2}-\gamma(p-1)u^{p-1}\triangle\log u<0
   \end{equation}
  holds everywhere. Actually, it can be inferred from \eqref{b.1} that there exist $\alpha\in(0,2/N)$ such that
   \begin{eqnarray}\nonumber\label{b.2}
    &-\Big(\frac{2}{N}-\alpha\Big)|\triangle\log u|^2-\frac{\alpha}{\gamma^2}\rho^2-\Bigg[\frac{\alpha(\gamma-\beta)^2}{\gamma^2}+(\gamma-\beta)(p-1)\Bigg]u^{2p-2}-\frac{\alpha(\gamma-1)^2}{\gamma^2}\Big|\frac{\nabla u}{u}\Big|^4&\\ \nonumber
    &-\gamma(p-1)u^{p-1}\triangle\log u-\Bigg[(p-1)(\gamma-2+\beta p)+\frac{2\alpha(\gamma-\beta)(\gamma-1)}{\gamma^2}\Bigg]u^{p-1}\Big|\frac{\nabla u}{u}\Big|^2&\\
    &-\frac{2\alpha(\gamma-\beta)}{\gamma^2}u^{p-1}\rho-\frac{2\alpha(\gamma-1)}{\gamma^2}\rho\Big|\frac{\nabla u}{u}\Big|^2\leq0.&
   \end{eqnarray}
Henceforth \eqref{e3.11} holds true. The condition \eqref{b.1} can be summarized as
    \begin{equation}\label{e3.12}
      \gamma>\max\{1,\beta,2-\beta p\},\ \ \ \frac{\gamma^2}{\gamma-\beta}(p-1)<\frac{8}{N}.
    \end{equation}
 Now, we divided \eqref{e3.12} into three cases:\\

 \noindent\textbf{Case 1:} If we take $\beta$ closing sufficiently to $0$, then for $1<p<1+\frac{4}{N}$, there exists $\gamma>2$ such that
   $$
    \frac{\gamma^2}{\gamma-\beta}(p-1)<\frac{8}{N}\Leftrightarrow\gamma(p-1)<\frac{8}{N},
   $$
 and thus the lemma holds true for $1<p<1+\frac{4}{N}$.\\

 \noindent\textbf{Case 2:} We can now assume $1+\frac{4}{N}\leq p<2$, and
  \begin{equation}\label{e3.13}
    \max\{1,\beta,2-\beta p\}=2-\beta p\Leftrightarrow\beta<1/p
  \end{equation}
 together with
   \begin{equation}\label{e3.14}
     2-\beta p>\gamma_0\equiv2\beta\Leftrightarrow\beta<\frac{2}{p+2}.
   \end{equation}
 Noting that the function $\frac{\gamma^2}{\gamma-\beta}$ attains its minimum $4\beta$ on $(\beta,+\infty)$ at the point $\gamma=\gamma_0\equiv2\beta$, upon the assumptions \eqref{e3.13} and \eqref{e3.14}, condition \eqref{e3.12} can be transformed into
    \begin{equation}\label{e3.15}
      \frac{(2-\beta p)^2}{2-\beta p-\beta}(p-1)<\frac{8}{N}\Leftrightarrow(p-1)(2-\beta p)^2-\frac{8(p+1)}{Np}(2-\beta p)+\frac{16}{Np}<0.
    \end{equation}
 Thus, it is not hard to see that for
    \begin{equation}\label{e3.16}
      \frac{N+2+\sqrt{N^2+12N+4}}{2N}\leq p<\min\Bigg\{2,\frac{N+2+\sqrt{N^2+8N}}{2(N-1)}\Bigg\},
    \end{equation}
 one can choose $\beta>0$ fulfilling \eqref{e3.13}-\eqref{e3.15}. Thus the lemma holds for \eqref{e3.16}.\\

 \noindent\textbf{Case 3:} Assuming that $p\geq2$, then the conditions \eqref{e3.13} and \eqref{e3.14} can be transformed into
   \begin{equation}\label{e3.17}
    0<\beta<\frac{1}{p}\Leftrightarrow1<2-\beta p<2.
   \end{equation}
 Therefore, the existence of $\beta>0$ satisfying \eqref{e3.15} and \eqref{e3.17} is equivalent to
   \begin{equation}\label{e3.18}
     \begin{cases}
       p<\frac{N+2+\sqrt{N^2+8N}}{2(N-1)},\\
       \frac{N+2+\sqrt{N^2+12N+4}}{2N}\leq p<\frac{N+4+\sqrt{N^2+24N+16}}{2N}.
     \end{cases}
   \end{equation}

 Summation above, because of
   $$
    \frac{N+2+\sqrt{N^2+12N+4}}{2N}<1+\frac{4}{N},
   $$
 we conclude that the lemma holds for $p$ satisfying \eqref{e3.9}. $\Box$\\

\vspace{40pt}

\section{Differential Harnack inequality of Li-Yau type}

Li-Yau type gradient estimation was first introduced by Li-Yau \cite{LY} for linear heat equation. This technic displayed a strengthened power in the past years and was later developed rapidly by vast literatures therein to nonlinear equations. See for examples \cite{CLL,D,HZ,L,LZ,MZ,MZS,W}. In this section, we will prove the following differential Harnack inequality of Li-Yau type:

\begin{prop}\label{p4.1} Let $(M,g)$ be a $N$ dimensional complete noncompact manifold with the metric $g$, whose Ricci curvature is bounded from below by $-Kg_{ij}$ for some nonnegative constant $K$. Assume that $u$ is a positive entire solution to \eqref{e1.4} and $p$ be an exponent satisfying \eqref{e3.9}, if $\beta,\gamma$ are parameters coming from Lemma \ref{l3.1}, we have
    \begin{equation}\label{e4.1}
       \rho\equiv\Bigg|\frac{\nabla u}{u}\Bigg|^2-\gamma(\log u)_t+\beta u^{p-1}\leq C_0K
    \end{equation}
  holds for some positive constant $C_0$ depending only on $N, p, \beta, \gamma$.
\end{prop}

To prove the proposition, we need the following two lemmas.

\begin{lemm}\label{l4.1}
   Under the assumption of Proposition \ref{p4.1}, if $\rho$ is positive at some point $z_0\equiv(x_0,t_0)\in M\times(-\infty,+\infty)$, then
     \begin{equation}\label{e4.2}
       |\triangle\log u|^2\geq\frac{1}{\gamma^2}\rho^2+\Bigg(\frac{\gamma-\beta}{\gamma}\Bigg)^2u^{2p-2}+\Bigg(\frac{\gamma-1}{\gamma}\Bigg)^2\Bigg|\frac{\nabla u}{u}\Bigg|^4
     \end{equation}
   holds at $z_0$.
\end{lemm}

\noindent\textbf{Proof.} By \eqref{e3.2}, we have
   \begin{eqnarray*}
     |\triangle\log u|^2&=&\Bigg|(\log u)_t-u^{p-1}-\Big|\frac{\nabla u}{u}\Big|^2\Bigg|^2\\
      &=&\Bigg|-\frac{1}{\gamma}\rho-\frac{\gamma-\beta}{\gamma}u^{p-1}-\frac{\gamma-1}{\gamma}\Big|\frac{\nabla u}{u}\Big|^2\Bigg|^2.
   \end{eqnarray*}
Since $\rho(z_0)>0$ and
   $$
    \gamma>1,\ \ \  \gamma>\beta,
   $$
inequality \eqref{e4.2} is clear true. $\Box$\\

A second lemma to be used is the following Laplace comparison property in differential geometry, see for example \cite{SY}:

\begin{lemm}\label{l4.2}
  Let $(M,g)$ be a Riemannian manifold satisfying $R_{ij}\geq-Kg_{ij}$ for some nonnegative constant $K$, then lying inside the cut-locus of a given point $o\in M$, the distance function $r(x)\equiv dist(x,o), x\in M$ satisfying $|\nabla r|_g=1$ and
   \begin{equation}\label{e4.3}
     r\triangle_gr\leq(N-1)\Bigg(1+\sqrt{\frac{K}{N-1}}r\Bigg).
   \end{equation}
\end{lemm}

\noindent\textbf{Remark.} Utilizing an approximation argument as in \cite{SY} or \cite{LJ,M}, we may assume that the distance function $r(x)$ from Lemma \ref{l4.2} is smooth.\\

Fixing a point $o\in M$, for any $R\geq1$, let's take a standard cut-off function $\overline{\eta}\in C^\infty_0({\mathbb{R}})$ satisfying that
   $$
    0\leq\overline{\eta}\leq1, \ \ \overline{\eta}(\tau)=\begin{cases}
      1, & \tau\in[-R^2/2,R^2/2],\\
      0, & \tau\not\in(-R^2,R^2)
    \end{cases}
   $$
and
   $$
     |\overline{\eta}'|\leq\frac{C_0}{R^2},\ \ \ |\overline{\eta}''|\leq\frac{C_0}{R^4}
   $$
for some universal constant $C_0>0$. Letting $\kappa$ be a large integer to be chosen later, we set $\eta(x,t)=\overline{\eta}^\kappa(r^2(x))\overline{\eta}^\kappa(t)$ for distant function $r(x)$ coming from Lemma \ref{l4.2}. Multiplying \eqref{e3.11} by $\eta$, we have
   \begin{eqnarray*}
     (\partial_t-\triangle)(\rho\eta)&=&\eta(\rho_t-\triangle\rho)+\rho(\eta_t-\triangle\eta)-2\nabla\rho\nabla\eta\\
      &\leq&\eta\Bigg\{-\alpha_4|\triangle\log u|^2+2\nabla\log u\nabla\rho-\alpha_5\rho^2-\alpha_6\Big|\frac{\nabla u}{u}\Big|^4+2K\Big|\frac{\nabla u}{u}\Big|^2\Bigg\}\\
      &&-2\nabla(\rho\eta)\frac{\nabla\eta}{\eta}+\rho\Bigg(\eta_t-\triangle\eta+2\frac{|\nabla\eta|^2}{\eta}\Bigg).
   \end{eqnarray*}
If $\rho$ is non-positive everywhere, Proposition \ref{p4.1} is clear true. Otherwise, evaluating $\rho\eta$ at its positive maximum point $z_0\equiv(x_0,t_0)\in B_{R^2}(o)\times(-R^2,R^2)$, we get
  \begin{equation}\label{e4.4}
    \alpha_5\rho^2\eta+\alpha_6\Bigg|\frac{\nabla u}{u}\Bigg|^4\eta\leq 2K\Bigg|\frac{\nabla u}{u}\Bigg|^2\eta-2\rho\nabla\log u\nabla\eta+\rho\Bigg(\eta_t-\triangle\eta+2\frac{|\nabla\eta|^2}{\eta}\Bigg)
  \end{equation}
at $z_0$, where
   $$
    0=2\nabla\log u\nabla(\rho\eta)=2\eta\nabla\log u\nabla\rho+2\rho\nabla\log u\nabla\eta
   $$
has been used at this point. By Lemma \ref{l4.2}, there hold
  \begin{eqnarray}\nonumber\label{e4.5}
   |\nabla\eta|&=&2\kappa\overline{\eta}^{\kappa-1}(r^2(x))\overline{\eta}^\kappa(t)\big|\overline{\eta}'(r^2(x))r\nabla r\big|\leq C\eta^{\frac{\kappa-1}{\kappa}}R^{-1},\\
   |\eta_t|&=&\kappa\overline{\eta}^\kappa(r^2(x))\overline{\eta}^{\kappa-1}(t)|\overline{\eta}'(t)|\leq C\eta^{\frac{\kappa-1}{\kappa}}R^{-2}
  \end{eqnarray}
and
   \begin{eqnarray}\nonumber\label{e4.6}
     |\triangle\eta|&=&2\kappa\overline{\eta}^{\kappa-1}(r^2(x))\overline{\eta}^\kappa(t)\Big|\overline{\eta}'(r^2(x))(r\triangle_gr+|\nabla r|^2)+2\overline{\eta}''(r^2(x))r^2|\nabla r|^2\Big|\\ \nonumber
     &&+4\kappa(\kappa-1)\overline{\eta}^{\kappa-2}(r^2(x))\overline{\eta}^\kappa(t)\Big|\overline{\eta}'(r^2(x))r\nabla r\Big|^2\\
     &\leq&C\Big(\eta^{\frac{\kappa-1}{\kappa}}R^{-1}+\eta^{\frac{\kappa-2}{\kappa}}R^{-2}\Big)
   \end{eqnarray}
for some positive constant $C$. Substituting into \eqref{e4.4} and using the Young's inequality, we obtain that
   \begin{eqnarray*}
     \alpha_5\rho^2\eta+\alpha_6\Bigg|\frac{\nabla u}{u}\Bigg|^4\eta\leq \varepsilon\Bigg(\rho^2+\Big|\frac{\nabla u}{u}\Big|^4\Bigg)\eta+ C_{0\varepsilon}K^2\eta+C_\varepsilon\Big(\eta^{\frac{\kappa-4}{\kappa}}R^{-4}+\eta^{\frac{\kappa-2}{\kappa}}R^{-2}\Big)
   \end{eqnarray*}
for any small constant $\varepsilon>0$, where $C_{0\varepsilon}$ and $C_\varepsilon$ are large constants independing and depending on $K$ respectively. Taking $\varepsilon$ small and then letting $R\to+\infty$, we conclude that $\rho\leq C_0K$ holds everywhere if $\kappa\geq4$ and the proof was done. $\Box$\\

\vspace{40pt}

\section{Non-existence of positive entire solution for nonlinear equation}

In this section, we prove the following non-existence result of positive entire solution to \eqref{e1.4}:

\begin{theo}\label{t5.1}
  Let $(M,g)$ be a $N$ dimensional complete noncompact manifold equipped with the metric $g$ and $u$ be a nonnegative entire solution to \eqref{e1.4}. Suppose that the exponent $p$ satisfies the subcritical condition \eqref{e3.9} and the Ricci curvature satisfies
   $$
    \begin{cases}
       R_{ij}\geq0, & \mbox{ for general manifold,}\\
       R_{ij}\geq -Kg_{ij}, & \mbox{ for suitable manifold}
    \end{cases}
   $$
for some positive constant $K$, then $u$ must be identical to zero.
\end{theo}

The first case for nonnegative Ricci curvature is easy. In fact, if $u$ is not identical to zero, by strong maximum principle, $u$ must be positive everywhere. Applying Proposition \ref{p4.1} for $K=0$, we have
   $$
    u_t\geq\frac{\beta}{\gamma}u^p, \ \ \forall (x,t)\in M\times(-\infty,+\infty).
   $$
Solving the above ordinary differential inequality, one can obtain that
   $$
    u(x,t)\geq\Bigg(u^{1-p}(x,0)-\frac{\beta}{\gamma}(p-1)t\Bigg)^{-\frac{1}{p-1}}, \ \ \forall t>0,
   $$
which is impossible since $u$ will blow up in finite time.

To prove the second case $R_{ij}\geq-Kg_{ij}$ for some positive constant $K$, we need the following lemma:

\begin{lemm}\label{l5.1}
  Let $(M,g)$ be a $N$ dimensional complete noncompact manifold equipped with the metric $g$ and $p$ satisfying the subcritical condition \eqref{e3.9}. Suppose that $R_{ij}\geq-Kg_{ij}$ holds for some positive constant $K$, then every positive entire solution to \eqref{e1.4} must be bounded from above by a positive constant
     $$
       C_{N,p,K}\equiv\Bigg(\frac{2C_0K}{\beta}\Bigg)^{\frac{1}{p-1}},
     $$
  where $\beta$ is a constant coming from Lemma \ref{l3.1}.
\end{lemm}

\noindent\textbf{Proof.} Applying Proposition \ref{p4.1} to positive $K$, one can obtain that
    \begin{equation}\label{e5.1}
     u_t\geq\frac{\beta}{\gamma}u^p-\frac{C_0K}{\gamma}u, \ \ \forall (x,t)\in M\times(-\infty,+\infty).
    \end{equation}
If the conclusion is not true, then there exists some point $(x_0,t_0)\in M\times(-\infty,+\infty)$, such that $u(x_0,t_0)> C_{N,p,K}$. Then, we claim that
   \begin{equation}\label{e5.2}
     u(x_0,t)> C_{N,p,K}, \ \ \forall t\geq t_0.
   \end{equation}
In fact, whenever $u(x_0,t)> C_{N,p,K}$, we have
   \begin{equation}\label{e5.3}
    u_t(x_0,t)>\frac{\beta}{2\gamma}u^p(x_0,t)>0.
   \end{equation}
So, \eqref{e5.2} keeps to hold for all $t\geq t_0$. Thus, back to turn, \eqref{e5.3} keeps to hold for all $t\geq t_0$. After solving the ordinary differential inequality \eqref{e5.3}, one can prove that
   $$
    u(x_0,t)\geq\Bigg(u^{1-p}(x_0,t_0)-\frac{\beta}{2\gamma}(p-1)(t-t_0)\Bigg)^{-\frac{1}{p-1}}, \ \ \forall t\geq t_0,
   $$
which is contradicting with $u$ is a entire solution. The proof was done. $\Box$\\

Now, one can complete the proof of Theorem \ref{t5.1}. Noting that when $M$ is a suitable manifold, the transformation
   $$
    u_k(x,t)=k^{\frac{2}{p-1}}u(kx,k^2t)
   $$
keeps the equation \eqref{e1.4} invariant for any given $k>0$. So, if there is a positive entire solution to \eqref{e1.4}, then the functions $u_k$ defined in above form an unbounded sequence of entire solutions of \eqref{e1.4}, contradicting with Lemma \ref{l5.1}. The conclusion of Theorem \ref{t5.1} was drawn. $\Box$\\

\vspace{40pt}

\section{Simple proofs to Yau's theorem and Gidas-Spruck's theorem}

In this final section, we prove firstly a Li-Yau type inequality for linear heat equation \eqref{e1.2}. In fact, a similar computation as \eqref{e3.2}-\eqref{e3.4} shows that
  \begin{eqnarray}\nonumber\label{e6.1}
    &&v_t-\triangle v=|\nabla v|^2,\\
    &&(\partial_t-\triangle)w=-2|\nabla^2v|^2-2Ric(\nabla v,\nabla v)+2\nabla v\nabla w,\\ \nonumber
    && (\partial_t-\triangle)(\log u)_t=2\nabla v\nabla(\log u)_t
  \end{eqnarray}
for positive solution $u$ to \eqref{e1.2}, where $v\equiv\log u, w\equiv|\nabla v|^2$. Letting $\gamma>1$ be a constant and setting $\rho\equiv w-\gamma(\log u)_t$, it is inferred from \eqref{e6.1} that
  \begin{equation}\label{e6.2}
    \rho_t-\triangle\rho=-2|\nabla^2v|^2-2Ric(\nabla v,\nabla v)+2\nabla v\nabla\rho.
  \end{equation}

 To show the desired differential Harnack inequality, we need Lemma \ref{l4.2} and the following variant version of Lemma \ref{l4.1}:\\

\begin{lemm}\label{l6.1}
  Let $u$ be a positive solution to \eqref{e1.2} and $v, w, \rho$ be defined as above. Then
    \begin{equation}\label{e6.3}
      N|\nabla^2v|^2\geq|\triangle v|^2\geq\frac{1}{\gamma^2}\rho+\Bigg(\frac{\gamma-1}{\gamma}\Bigg)^2w^2
    \end{equation}
  holds at the place where $\rho$ is positive.
\end{lemm}

\noindent\textbf{Proof.} The first inequality comes from Cauchy's inequality. The second inequality is also element since
  \begin{eqnarray*}
   |\triangle v|^2&=&\big|(\log u)_t-w\big|^2=\Bigg|-\frac{1}{\gamma}\rho-\frac{\gamma-1}{\gamma}w\Bigg|^2\\
    &\geq&\frac{1}{\gamma^2}\rho^2+\Bigg(\frac{\gamma-1}{\gamma}\Bigg)^2w^2,
  \end{eqnarray*}
as long as $\rho$ is positive. $\Box$\\

Now, for any given $R>1$, we take a cut-off function $\eta$ as in Section 4. Applying \eqref{e4.5} and evaluating $\rho\eta$ at its possible positive maximum point $z_0$, one can verify that
  \begin{eqnarray*}
    0&\leq&(\partial_t-\triangle)(\rho\eta)=\eta(\rho_t-\triangle\rho)+\rho(\eta_t-\triangle\eta)-2\nabla\rho\nabla\eta\\
     &\leq&\eta\Bigg\{-\frac{2}{N\gamma^2}\rho^2-\frac{2}{N}\Bigg(\frac{\gamma-1}{\gamma}\Bigg)^2w^2+2Kw+2\nabla v\nabla\rho\Bigg\}+\rho\Bigg(\eta_t-\triangle\eta+2\frac{|\nabla\eta|^2}{\eta}\Bigg)\\
     &\leq&-\frac{2}{N\gamma^2}\rho^2\eta-\frac{1}{N}\Bigg(\frac{\gamma-1}{\gamma}\Bigg)^2w^2\eta+CK^2\eta+2\rho\nabla v\nabla\eta+\rho\Bigg(\eta_t-\triangle\eta+2\frac{|\nabla\eta|^2}{\eta}\Bigg)\\
     &\leq&-\frac{1}{N\gamma^2}\rho^2\eta-\frac{1}{2N}\Bigg(\frac{\gamma-1}{\gamma}\Bigg)^2w^2\eta+CK^2\eta+C\Big(\eta^{\frac{\kappa-4}{\kappa}}R^{-4}+\eta^{\frac{\kappa-2}{\kappa}}R^{-2}\Big),
  \end{eqnarray*}
where $C$ is a positive constant depending only on $N,\gamma$ and
   \begin{eqnarray*}
     0&=&-2\nabla(\rho\eta)\frac{\nabla\eta}{\eta}=-2\nabla\rho\nabla\eta-2\rho\frac{|\nabla\eta|^2}{\eta},\\
     0&=&2\nabla v\nabla(\rho\eta)=2\eta\nabla v\nabla\rho+2\rho\nabla v\nabla\eta
   \end{eqnarray*}
has been used. So, we obtain that
   \begin{eqnarray*}
     \rho^2\eta^2&\leq&\rho^2\eta\big|_{z_0}\leq C(K^2+R^{-4}+R^{-2}).
   \end{eqnarray*}
Sending $R\to+\infty$, we get the following differential Harnack inequality:\\

\begin{theo}\label{t6.1}
  Let $(M,g)$ be a $N$ dimensional complete noncompact manifold with the metric $g$. Suppose that the Ricci curvature satisfies $R_{ij}\geq-Kg_{ij}$ for some nonnegative constant $K$, then every positive entire solution $u$ of \eqref{e1.2} satisfies that
    \begin{equation}\label{e6.4}
      \Bigg|\frac{\nabla u}{u}\Bigg|^2-\gamma(\log u)_t\leq C_{N,\gamma}K
    \end{equation}
  for any $\gamma>1$, where $C_{N,\gamma}$ is a positive constant depending only on $N, \gamma$.
\end{theo}

Now, Yau's theorem and Gidas-Spruck's theorem are direct consequences of Theorem \ref{t6.1} and Proposition \ref{p4.1} respectively:\\

\begin{theo}\label{t6.2} Let $(M,g)$ be a $N$ dimensional complete noncompact manifold with the metric $g$. We have the following alternatives:

 (1) Suppose that $1<p<p_*(N)$ and the Ricci curvature satisfying
   $$
    \begin{cases}
       R_{ij}\geq0, & M\ \mbox{is a general manifold,}\\
       R_{ij}\geq -Kg_{ij}, & M\ \mbox{is a suitable manifold}
    \end{cases}
   $$
for positive constant $K$. Then there is no positive entire solution to \eqref{e1.3}, and

 (2) assuming only $R_{ij}\geq0$, every nonnegative entire solution to \eqref{e1.1} must be a constant.\\
\end{theo}

\noindent\textbf{Remark.} Part (1) is only a special case of the theorem of Gidas-Spruck \cite{GS}. We present here using a different proof.\\

\noindent\textbf{Proof.} Part (1) is a corollary of Theorem \ref{t5.1} since every solution to \eqref{e1.3} is a stationary solution to \eqref{e1.4}. Part (2) follows from \eqref{e6.4}. In fact, since $u_t\equiv0$ and $K=0$, it is inferred from \eqref{e6.4} that
   $$
    \nabla u(x)\equiv 0,\ \ \forall x\in M.
   $$
 Thus the conclusion is drawn. $\Box$\\

\vspace{40pt}

\section*{Acknowledgments}

The author (SZ) would like to express his deepest gratitude to Professors Xi-Ping Zhu, Kai-Seng Chou, Xu-Jia Wang and Neil Trudinger for their constant encouragements and warm-hearted helps. Special thanks were also owed to Professor Jun-Cheng Wei and Professor Ke-Lei Wang for valuable suggestions and helpful conversations. This paper was also dedicated to the memory of Professor Dong-Gao Deng.\\

\end{document}